\documentclass[12pt]{amsart}
\usepackage{amsmath}
\usepackage{amssymb}
\newtheorem{theorem}{Theorem}[section]
\newtheorem{lemma}[theorem]{Lemma}

\newtheorem{corollary}[theorem]{Corollary}
\newtheorem{claim}[theorem]{Claim}

\theoremstyle{definition}
\newtheorem{definition}[theorem]{Definition}
\newtheorem{definitions}[theorem]{Definitions}

\newtheorem{definitions and remarks}[theorem]{Definitions and Remarks}
\newtheorem{question}[theorem]{Question}

\theoremstyle{remark}
\newtheorem{remark}[theorem]{Remark}

\newtheorem{exercises}[theorem]{Exercises}

\numberwithin{equation}{section}

\newcommand{\exc}{\mathrm{exc}}
\newcommand{\inv}{\mathrm{inv}}
\newcommand{\oinv}{\overline{\inv}}
\newcommand{\omu}{\overline{\mu}}
\newcommand{\oJ}{\overline{J}}

\newcommand{\Sing}{\mathrm{Sing}\,}
\newcommand{\supp}{\mathrm{supp}\,}
\newcommand{\cosupp}{\mathrm{cosupp}\,}

\newcommand{\ord}{\mathrm{ord}}
\newcommand{\order}{\mathrm{order}}

\newcommand{\length}{\mathrm{length}\,}

\newcommand{\Der}{\mathrm{Der}}
\newcommand{\lcm}{\mathrm{lcm}}
\newcommand{\Sub}{\mathrm{Sub}}

\newcommand{\al}{{\alpha}}

\newcommand{\de}{{\delta}}

\newcommand{\g}{{\gamma}}
\newcommand{\G}{{\Gamma}}

\newcommand{\La}{{\Lambda}}

\newcommand{\p}{{\partial}}
\newcommand{\s}{{\sigma}}

\newcommand{\vp}{{\varphi}}
\newcommand{\io}{{\iota}}

\newcommand{\IN}{{\mathbb N}}
\newcommand{\IQ}{{\mathbb Q}}
\newcommand{\IA}{{\mathbb A}}

\newcommand{\cC}{{\mathcal C}}
\newcommand{\cD}{{\mathcal D}}
\newcommand{\cE}{{\mathcal E}}

\newcommand{\cG}{{\mathcal G}}
\newcommand{\cH}{{\mathcal H}}
\newcommand{\cI}{{\mathcal I}}
\newcommand{\cJ}{{\mathcal J}}
\newcommand{\cM}{{\mathcal M}}
\newcommand{\cN}{{\mathcal N}}
\newcommand{\cO}{{\mathcal O}}

\newcommand{\cW}{{\mathcal W}}

\newcommand{\um}{{\mathfrak m}}

\newcommand{\ucG}{\underline{\cG}}
\newcommand{\ucH}{\underline{\cH}}
\newcommand{\ucI}{\underline{\cI}}
\newcommand{\ucJ}{\underline{\cJ}}
\newcommand{\ucC}{\underline{\cC}}
\newcommand{\ucD}{\underline{\cD}}
\newcommand{\ucE}{\underline{\cE}}
\newcommand{\ucM}{\underline{\cM}}
\newcommand{\ucN}{\underline{\cN}}

\newcommand{\ucW}{\underline{\cW}}

\newcommand{\uX}{\underline{X}}

\begin{document}

\title[Functoriality in resolution of singularities]
{Functoriality in resolution of singularities}

\dedicatory{To Heisuke Hironaka for his 77th birthday,\\
in celebration of \emph{Kiju} --- joy and long life!}

\author{Edward Bierstone}
\address{Department of Mathematics, University of Toronto, 40 St. George Street,
Toronto, Ontario, Canada M5S 2E4}
\email{bierston@math.toronto.edu}
\thanks{The authors' research was supported in part by NSERC 
grants OGP0009070 and OGP0008949.}

\author{Pierre D. Milman}
\address{Department of Mathematics, University of Toronto, 40 St. George Street,
Toronto, Ontario, Canada M5S 2E4}
\email{milman@math.toronto.edu}



%\begin{abstract}
%\end{abstract}

\maketitle
\setcounter{tocdepth}{1}
\tableofcontents

\section{Introduction}
Hironaka's theorem on resolution of singularities in
characteristic zero \cite[1964]{Hann} is one of the monuments of twentieth
century mathematics. Canonical or functorial versions of Hironaka's theorem
were established 
more than fifteen years ago, by Villamayor \cite{V1, V2} and by the authors
\cite{BMmega, BMinv}. Two new treatments of canonical resolution
of singularities have appeared just in the last couple of years, by 
W{\l}odarcyzk \cite{W} and Kollar \cite{Ko}. These papers reflect a
continuing interest in a better understanding of desingularization,
in large part with the goal of discovering techniques that extend to positive
characteristic (cf. \cite{Hinf, Hprog, Ka}). They also show that certain 
aspects have remained mysterious, even in characteristic zero. Our goal here
is to clarify these issues (and raise some questions).

The various proofs of desingularization have a lot in common, but also
important differences that are related to
the notion of {\it functoriality} or {\it canonicity} involved. 
They all involve reducing the problem to
``canonical desingularization'' of a collection of local resolution data,
originating in Hironaka's idea of an {\it idealistic exponent}, and called
a {\it presentation} \cite{BMinv}, {\it basic object} \cite{V2} or 
{\it marked ideal} \cite{W} (with variations in the meaning of these
objects). In all cases, canonical desingularization is proved by induction
on dimension. Functoriality (apart from its intrinsic interest) plays two
important roles: (i) in the statements of the theorems that can be proved
as consequences of canonical desingularization of a marked ideal; (ii) in the
proof itself --- it is sometimes easier to prove a stronger theorem by
induction, because we can make a stronger inductive assumption. (We use
W{\l}odarcyzk's terminology ``marked ideal'', though 
our notion is a little more general. See \S1.1.)

We begin with the theorem of \emph{canonical resolution of singularities} in
characteristic zero (not the most general statement, but formulated in a way 
that is easy to state and includes the most useful conditions). 
For simplicity, we restrict our attention to algebraic varieties (or reduced
separated schemes of finite type). (There
are more general theorems for nonreduced schemes; see \cite{Hann}, 
\cite[Sect.\,11]{BMinv}. Moreover, all results below apply to analytic spaces,
where the sequences of blowings-up are finite over any relatively compact
open set.)

\begin{theorem}
Given an algebraic variety $X$ over a field of characteristic zero, 
there is finite sequence of blowings-up
$\s_j$ with smooth centres,
$$
X = X_0 \stackrel{\s_1}{\longleftarrow} X_1 \longleftarrow \cdots
\stackrel{\s_{t}}{\longleftarrow} X_{t}\,,
$$
such that:
\begin{enumerate}
\item 
$X_t$ is smooth and the exceptional divisor $E_t$ in $X_t$ has only 
normal crossings.
\item
The morphism $\s$ given by the composite of the $\s_j$ is an isomorphism 
over $X \setminus \Sing X$.
\item
The resolution morphism $\s_X=\s: X_t \to X$ (or 
the entire sequence of blowings-up $\s_j$) can be associated to $X$
in a way that is functorial (at least with 
respect to \'etale or smooth morphisms, and 
field extensions).
\end{enumerate}
\end{theorem}

This theorem can be proved with following stronger version of the
condition (2):
For each $j$, let $C_j \subset X_j$ denote the centre of the blowing-up
$\s_{j+1}: X_{j+1} \to X_j$, and let $E_j$ denote the exceptional divisor
of $\s_1\circ\cdots\circ\s_j$. Then:
\begin{enumerate}
\item[($2'$)] {\it For each $j$, either $C_j \subset \Sing X_j$
or $X_j$ is smooth and $C_j \subset \supp E_j$.}
\end{enumerate}

Theorem 1.1 is proved in \cite{BMinv}, in \cite{Hann} without
the functoriality condition (3), and in \cite{V2} assuming Hironaka's
reduction to the case of an idealistic exponent \cite{Hid}.
A weaker version of the theorem
is proved in \cite{EVnew, W, Ko},
where the blowings-up $\s_j$ do not necessarily have smooth centres.
(We recall that every birational morphism
of quasi-projective varieties is a blowing-up with centre that is not
necessarily smooth \cite[Thm.\,7.17]{Hart}.) We do not know of a proof
that provides the condition of smooth centres without also giving ($2'$).
(See \S1.3 below.)

There is no version of Theorem 1.1 that is functorial with
respect to all morphisms. Consider the morphism $\vp$ from
$\IA^2$ to the quadratic cone $X := \{uv-w^2=0\} \subset \IA^3$ given
by $\vp(x,y)=(x^2,y^2,xy)$ \cite[Ex.\,3.4]{Ko}.
Then $\s_{\IA^2}$ is the identity morphism
(by functoriality with respect to translations), so functoriality with
respect to $\vp$ would imply that $\vp$ lifts to the blowing-up $\s:
X' \to X$ over the origin (for example, by the universal mapping property
of blowing up). But this is clearly false.

\subsection{Reduction to desingularization of a marked ideal}
For the sake of an inductive proof of Theorem 1.1,
it is natural to deduce the result from an {\it embedded} version:
We assume that $X \subset M$, where $M$ is smooth. Then the successive
blowings-up of $X$ are given by  restricting a sequence of blowings-up
of $M$ to the corresponding strict transforms of $X$. 

All currently available proofs of canonical embedded resolution of singularities
follow Hironaka's idea of reducing the result to ``functorial desingularization''
of an ideal with associated multiplicity \cite{Hid}:

A {\it marked ideal} is a quintuple $\ucI = (M,N,E,\cI,d)$, where $N$ is a smooth
subvariety of $M$, $E$ is a normal crossings divisor, $\cI$ is a
coherent ideal (sheaf) on $N$, and $d \in \IN$. 
$\ucI$ encodes the 
data we have at each step of resolution: $N$ will be a
{\it maximal contact} subvariety and $E$ an exceptional divisor.
A marked ideal has a 
transformation law that is simpler than strict transform:
Say that a blowing-up
$\s: M'\to M$ is {\it admissible} if its centre $C \subset \cosupp \ucI
:= \{x \in N: \order_x\cI \geq d\}$, and $C,\,E$ have only normal crossings.
The {\it (weak) transform} $\cI'$ of $\cI$ by $\s$ 
is generated by $y_\exc^{-d} f \circ \s$,
for all $f \in \cI$, where $y_\exc$ is a local generator of the ideal
of $\s^{-1}(C)$. We define a {\it resolution of
singularities} of $\ucI$ as a finite sequence of admissible blowings-up,
after which the cosupport is empty. (See Section 2 for detailed definitions.)

The definition of a marked ideal above is a slight generalization of that of
\cite{W} and corresponds to a \emph{presentation} in \cite{BMinv}.
W{\l}odarczyk's marked ideal
would be given by $(N, E|_N, \cI,d)$.

We prove desingularization of a marked ideal by induction
on $\dim N$: We reduce the dimension by passing to an {\it equivalent
coefficient ideal} 
on a smooth {\it maximal contact hypersurface}. (See Defns. 2.5
and Section 4.) Maximal contact exists
locally in characteristic zero, but is not unique.
Coefficient ideals are defined using ideals of derivatives of $\cI$
(The notion goes back to \cite{Hid, giraud}, and is used in \cite{V2}). 
An idea in \cite{BMinv, BMda1} was to use 
derivatives preserving the exceptional
divisor $E$ (i.e., derivatives logarithmic with respect to the exceptional
variables). The idea is natural to the strong version of
functoriality in Theorem 1.3 below (see \S3.4), and 
is crucial to Kawanoue's proposal for 
desingularization in positive characteristic \cite{Ka}. The use of derivatives
preserving $E$ in resolution of singularities occurs already in \cite{C}.

Given $X \subset M$, let $\cI_X \subset \cO_M$ denote the ideal of $X$.
Desingularization of the marked ideal
$\ucI = (M, M, \emptyset, \cI_X, 1)$
induces embedded resolution of 
singularities of $X$
in the case that $X \subset M$ is a hypersurface (i.e., of codimension $1$),
but only a weaker theorem in general (cf. \cite{EVnew, W}):
Consider the sequence of blowings-up $\{\s_j\}$ of $M$ that resolves
the singularities of $\ucI$ (and let $\ucI_j$ denote the corresponding 
transforms of $\ucI$). For some $r$, the centre $C_r$ of 
$\s_{r+1}$
coincides with $\cosupp \ucI_r$. Then $X_r$ is smooth
and has only normal crossings with $E_r$ (by Defns. 2.2). 
We get conditions (1)--(3) of
Theorem 1.1, but the intersections $C_j \cap X_j$ are not necessarily
smooth; i.e., the induced sequence of blowings-up of $X$ does not necessarily
have smooth centres (nor do we necessarily get
$C_j \cap X_j \subset \Sing X_j$).
See Example 8.2. To get
Theorem 1.1, we use a stronger idea of functoriality, and {\it equisingularity}
of the centres of blowing up with respect to the {\it Hilbert-Samuel function}.
We know of no other way, contrary to
Kollar's suggestion \cite[p.2]{Ko}. (See \S1.3.)

\subsection{Test sequences and functorial desingularization}
Hironaka \cite{Hid} already proposed to deal 
with the preceding issue by 
using blowings-up not only as building blocks in resolution of
singularities, but also to test for equisingularity of marked ideals.
Define a {\it weak test sequence} for $\ucI$ 
as a sequence of morphisms, each either an admissible blowing-up
or a projection $M \times \IA^1 \to M$. Say that
two marked ideals on $M$ are {\it weak equivalent} if they have
the same weak test sequences (Defns. 2.4).  

Functoriality in a general sense should mean that ``equivalent'' marked ideals
undergo the same resolution process.
In \cite[\S9]{Hid}, Hironaka asserts (without proof) that a marked ideal 
can be desingularized
functorially, not only with respect to smooth morphisms,
but also with respect to weak equivalence classes (i.e., weak equivalent
marked ideals have the same resolution sequences). 

\begin{question} 
Is this true?
\end{question}

The question is interesting partly in view of the role played by
weak equivalence in \cite{Hinf}, which develops techniques
for use in positive characteristic. We do not even know whether the
algorithms of \cite{BMinv, EVnew, W, Ko} are functorial in the sense of
Hironaka's assertion.

In \cite{BMinv, BMda1}, we realize Hironaka's philosophy by
enlarging the class of test transformations to include
{\it exceptional blowing-ups} (centre $=$
an intersection of two components of the exceptional divisor). Say that
two marked ideals on $M$ are {\it equivalent} if they have the same
sequences of test transformations. (See Defns. 2.5.)

\begin{theorem}[Functorial resolution of singularities of a marked ideal]
Let $\ucI = (M,N,E,\cI,d)$ denote a marked ideal. Set $\ucI_0 := \ucI$, $\ucI_0 =
(M_0,N_0,E_0,\cI_0,d_0)$. Then there is a finite sequence of blowings-up
with smooth centres,
\begin{equation}
M = M_0 \stackrel{\s_1}{\longleftarrow} M_1 \longleftarrow \cdots
\stackrel{\s_{t}}{\longleftarrow} M_{t}\,,
\end{equation}
where $\s_1$ is admissible for $\ucI_0$ and each 
successive $\s_{j+1}$ is admissible for the transform 
$\ucI_j$ of $\ucI_{j-1}$, such that
$$
\cosupp \ucI_t = \emptyset\,.
$$
Moroever, a resolution sequence (1.1) can be associated to every marked ideal
$\ucI$ in a way that is functorial with respect to equivalence classes
of marked ideals and smooth morphisms.
\end{theorem}

The \emph{functoriality} assertion here means that, if $\ucI = (M,N,E,\cI,d)$,
$\ucI_1 = (M_1,N_1,E_1,\cI_1,d_1)$ are marked ideals and $\vp: M_1 \to M$ 
is a smooth morphism
such that $\ucI_1$ is equivalent to $\vp^*(\ucI)$, then $\vp$ lifts to smooth
morphisms throughout the resolution towers for $\ucI_1$, $\ucI$. Functoriality
with respect to smooth morphisms in \cite{W, Ko} is the weaker assertion
that $\vp$ lifts to smooth morphisms throughout the resolution towers for
$\vp^*(\ucI)$ and $\ucI$. When we say that $\vp$ lifts to smooth morphisms
throughout the resolution sequences, we allow the possibility of inserting
trivial blowings-up (i.e., identity morphisms) in the resolution tower for
$\ucI_1$. (For example, if $\vp: U \hookrightarrow M$ is an embedding of
an open subvariety $U$, then the centre of a given blowing-up in the resolution
sequence for $\ucI$ may have no points over $U$.) This issue is treated in
Section 7.

We give a complete proof of Theorem 1.3 in Sections 5--7. The proof
is essentially a repetition of that in \cite{BMinv}, but we use
the very clear inductive framework of \cite{W}. Part
of our purpose is to clarify the relationship between various proofs,
particularly regarding the role of functoriality 
and the different notions of derivative ideal 
(Sect.\,5). We also want to isolate the precise role of a
desingularization invariant (Sect.\,7) and to show that the idea of
equivalence makes it possible to compare the actual algorithms for choosing
the centres of blowing up. In particular, we have the following corollary
of the proof of Theorem 1.3.

\begin{corollary}
The algorithm for resolution of singularities of marked ideals of 
W{\l}odarczyk \cite{W} or Kollar \cite{Ko} coincides
with that of \cite{BMinv}.
\end{corollary}

This is proved in \S8.1. We know of no way to prove it other
than using the strong notion of functoriality in Theorem 1.3. Our proof
of Theorem 1.3 shows that the association of a resolution sequence to
a marked ideal is also functorial with respect to field extensions and
to closed embeddings $\io: M \hookrightarrow M'$ in the case that $d=1$
(as in \cite{W, Ko}), but we do not explicitly discuss these issues.

Theorem 1.3 is proved by induction on $\dim \ucI := \dim N$. 
The inductive step breaks up into two distinct steps (Section 5):
\smallskip

\noindent
{\bf I.\ }
Functorial desingularization of a marked ideal of dimension $n-1$ implies 
functorial desingularization
of a marked ideal of \emph{maximal order}, of dimension $n$. 
(See Defns. 2.1.)
\smallskip

\noindent
{\bf II.\ }
Functorial desingularization of a marked ideal of maximal order of 
dimension $n$ implies functorial
desingularization of a general marked ideal of dimension $n$.
\smallskip

Step I involves passing to an
equivalent coefficient ideal on a local maximal
contact hypersurface. 
We have to show that the centres of blowing
up given by canonical desingularization of the coefficient ideals (which
holds by induction) in overlapping coordinate charts, glue together to define
a global centre of blowing up for the original marked ideal. Our proof 
and those \cite{W, Ko} treat this problem in rather different ways:
In our induction, the coefficient ideals in overlapping charts are equivalent,
so that Step I follows from the inductive assumption of strong
functorial desingularization in lower dimension -- this is the reason for
our notion of equivalence. But Step II requires some technical work. On the
other hand, in \cite{W, Ko}, most of the work is required for Step I (passage
to a ``homogenized ideal'', to reduce to a situation where maximal contact
is uniquely determined up to an \'etale automorphism), while II is simpler.

\subsection{Use of the Hilbert-Samuel function}
We can deduce Theorem 1.1 (including ($2'$)) from Theorem 1.3,
following Hironaka's philosophy in \cite{Hid}.
We use the \emph{Hilbert-Samuel function}
$$
H_{X,x}(k) := \length \frac{\cO_{X,x}}{\um_{X,x}^{k+1}}\,,\quad k \in \IN\,.
$$
The Hilbert-Samuel function $H_{X,\cdot}: X \to \IN^{\IN}$
has the following basic properties, established
by Bennett \cite{Be} (see \cite{BMjams, BMinv} for simple proofs): 
(1) $H_{X,\cdot}$ distinguishes smooth and singular points; (2) $H_{X,\cdot}$ is
Zariski upper-semicontinuous; (3) $H_{X,\cdot}$ is \emph{infinitesimally 
upper-semicontinuous} (i.e., $H_{X,\cdot}$ cannot increase after blowing-up
with centre on which it is constant); (4) Any decreasing sequence in the value
set of the Hilbert-Samuel function stabilizes.

The Hilbert-Samuel function $H_{X,\cdot}$ is a local invariant that 
plays the same role with respect to strict transform of a variety $X$
as the order with respect to (weak) transform of a (marked) ideal. More precisely,
for all $a \in X$, there is a local embedding $X \hookrightarrow M$ at $a$
and a marked ideal $\ucI = (M,N,\emptyset,\cI,d)$ which has the same
\emph{test sequences} as $\uX := (M,\emptyset,X,H)$, where $H = H_{X,a}$. (We
define a \emph{test sequence} for $\uX = (M,E,X,H)$ 
by analogy with that for a marked
ideal (Defns. 2.5), but where a blowing-up $\s: M' \to M$ with smooth centre $C$
is \emph{admissible} if $C \subset \supp \uX := \{x\in X: H_{X,x} \geq H\}$ and
$C,\,E$ have only normal crossings, and where $X$ transforms by \emph{strict
transform}.) We call $\ucI$ a \emph{presentation} of $H_{X,\cdot}$ at $a$.
If $a$ is a maximum point of $H_{X,\cdot}$, then the maximal value of
$H_{X,\cdot}$ decreases after desingularization of $\ucI$, and Theorem 1.1
follows (using the basic properties of $H_{X,\cdot}$ above and) using
functoriality in the same way as in Step I of the proof of Theorem 1.3.
(See \cite[Ch.\,III]{BMinv}.)

We are interested in equisingularity with respect to the Hilbert-Samuel 
function as a {\it necessary}
rather than only sufficient condition on the centres of
blowing up:

\begin{question}
Does any local invariant $\io_{X,\cdot}$ 
which satisfies the basic properties (1)--(4) above and which admits a
presentation determine the Hilbert-Samuel function $H_{X,\cdot}$?
For example, is the stratification of $X$
by the values of $\io_{X,\cdot}$ a refinement of that by $H_{X,\cdot}$?
\end{question}

\section{Marked ideals and test transformations}
In this section we give detailed definitions of the notions concerning
marked ideals, test transormations and equivalence that have been introduced
informally in Section 1.

\begin{definitions}
A \emph{marked ideal} $\ucI$ is a quintuple $\ucI = (M,N,E,\cI,d)$, where:
\begin{enumerate}
\item
$M$, $N$ are smooth, and $N \subset M$ is a closed subvariety;
\item
$E = \sum_{i=1}^s H_i$ is a simple normal crossings divisor on $M$ which is
tranverse to $N$ and \emph{ordered}. (The $H_i$ are smooth hypersurfaces
in $M$, not necessarily irreducible, with ordered index set as indicated.)
\item $\cI$ is an ideal (i.e., a coherent ideal sheaf) in $\cO_N$.
\item $d \in \IN$.
\end{enumerate}

We define the \emph{cosupport} $\cosupp \ucI$ of a marked ideal 
$\ucI = (M,N,E,\cI,d)$ as
$$
\cosupp \ucI := \{x \in N:\, \ord_x\cI \geq d\}\,,
$$
where $\ord_x\cI$ denotes the largest nonnegative integer $q$ such that
$\cI_x \subset \um^q_{N,x}$, with $\um_{N,x}$ the maximal ideal of $\cO_{N,x}$.
\smallskip

We say that a marked ideal $\ucI = (M,N,E,\cI,d)$ is of \emph{maximal
order} if $d = \max\{\ord_x \cI: x \in \cosupp \ucI\}$.
\end{definitions}

\begin{definitions}
Let $\ucI = (M,N,E,\cI,d)$ denote a marked ideal. A blowing-up
$\s: M' \to M$ (with smooth centre $C$) is \emph{admissible} for $\ucI$ if
$C \subset \cosupp \ucI$, and
$C$, $E$ have only normal crossings.
\smallskip

We define the \emph{transform} of $\ucI$ by an admissible blowing-up
$\s: M' \to M$ as the marked ideal $\ucI' = (M',N',E',\cI',d'=d)$,
where:
\begin{enumerate}
\item
$N'$ is the strict transform of $N$ by $\s$.
\item
$E' = \sum_{i=1}^{s+1} H_i'$, where $H_i'$ denotes the strict transform
of $H_i$, for each $i=1,\ldots,s$, and $H'_{s+1} := \s^{-1}(C)$ (the
exceptional divisor of $\s$, introduced as the last member of $E'$).
\item
$\cI' := \cI_{\s^{-1}(C)}^{-d}\cdot \s^*(\cI)$ (where $\cI_{\s^{-1}(C)}
\subset \cO_{N'}$ denotes the ideal of $\s^{-1}(C)$). I.e., $\cI'$ is
generated locally by $y_{\exc}^{-d} f\circ\s$, for all local generators
$f$ of $\cI$, where $y_{\exc}$ denotes a local generator of $\cI_{\s^{-1}(C)}$.
\end{enumerate}
\end{definitions}

In this definition, note that $\cI'$ is divisible by $\cI_{\s^{-1}(C)}^d$
and $E'$ is a normal crossings divisor transverse to $N'$, because $\s$
is admissible.

\begin{definition}
Let $\ucI = (M,N,E,\cI,d)$ denote a marked ideal, and let $\vp: M_1 \to M$
be a smooth morphism. We define the \emph{pull-back} $\vp^*\ucI$ of $\ucI$
as the marked ideal $\vp^*\ucI := (M_1, \vp^{-1}(N), \vp^{-1}(E), \vp^*(\cI), d)$.
\end{definition}

\begin{definitions} Let $\ucI = (M,N,E,\cI,d)$ denote a marked ideal as above.
Consider the product of $M$ with a line, $M' := M \times \IA^1$, and let
$\pi:\, M' \to M$ denote the projection morphism. We define the \emph{transform}
$\ucI'$ of $\ucI$ by $\pi$ as the marked ideal $\ucI' = (M',N',E',\cI',d'=d)$,
where $N' := \pi^{-1}(N)$, $\cI' := \pi^*(\cI)$, and $E' = 
\sum_{i=1}^{s+1} H_i'$, where $H_i' := \pi^{-1}(H_i)$, for each $i = 1,\ldots,s$,
and $H'_{s+1}$ denotes the \emph{horizontal divisor} $D := M \times \{0\}$
(included as the last member of $E'$). (Of course, $\pi$ is a smooth morphism,
but $\ucI'$ differs from $\pi^*(\ucI)$ in the term $E'$.)
\smallskip

We define a \emph{weak test sequence} for $\ucI_0 = \ucI$ as a sequence of morphisms
\begin{equation}
M = M_0 \stackrel{\s_1}{\longleftarrow} M_1 \longleftarrow \cdots
\stackrel{\s_{t}}{\longleftarrow} M_{t}\,,
\end{equation}
where each successive morphism $\s_{j+1}$ is either an admissible blowing-up
(for the transform $\ucI_j$ of $\ucI_{j-1}$), or the projection from the
product with a line.
\smallskip

We say that two marked ideals $\ucI$ and $\ucI_1 = (M_1 = M, N_1, E_1 = E,
\cI_1, d_1)$ (with the same ambient variety $M$ and the same normal crossings
divisor $E$) are \emph{weak-equivalent} if they have the same weak test sequences;
i.e., every weak test sequence for one is a weak test sequence for the other.
\end{definitions}

\begin{definitions}
Let $\ucI = (M,N,E,\cI,d)$ be a marked ideal. A blowing-up $\s: M' \to M$
will be called an \emph{exceptional blowing-up} for $\ucI$ if its centre
$C = H_i \cap H_j$, where $H_i$, $H_j$ are distinct hypersurfaces in $E$.
\smallskip

We define the \emph{transform} $\ucI'$ of $\ucI$ by an exceptional blowing-up
$\s: M' \to M$ (with centre $C$, say) as the marked ideal $\ucI'
= (M',N',E',\cI',d')$, where $N'$ the strict transform of $N$ ($N'=\s^{-1}(N)$,
in this case), $E'$ is defined as in Definitions 1.3 above, $\cI' = 
\s^*(\cI)$, and $d'=d$.
\smallskip

We define a \emph{test sequence} for $\ucI_0 = \ucI$ as a sequence of morphisms
(2.1), where each successive $\s_{j+1}$ is either an admissible blowing-up,
the projection from a product with a line, or an exceptional blowing-up. (Morphisms
of these three kinds will be called \emph{test transformations} or 
\emph{test morphisms}.)
\smallskip

We say that two marked ideals $\ucI$ and $\ucI_1$ (with the same ambient 
variety $M$ and the same normal crossings divisor $E$) are \emph{equivalent}
if they have the same test sequences. 
\end{definitions}

\section{Differential calculus}
Let $N$ denote a smooth variety. If $(x_1,\ldots,x_n)$ are local coordinates
at a point of $N$, then the partial derivatives $\p/\p x_1,\ldots, \p/\p x_n$
are well-defined as local generators of the sheaf of derivations 
$\Der_N$ of $\cO_N$

\subsection{Transformation of differential operators by blowing up}
Let $\s: N' \to N$ denote the blowing-up with centre a smooth subvariety
$C$ of $N$. Choose a local coordinate chart $U$ of $N$ with coordinates
$(x_1,\ldots,x_n)$ in which $C = \{x_r = \cdots = x_n =0\}$. Then $\s^{-1}(C)$
is covered by coordinate charts $U_{x_r},\ldots, U_{x_n}$, where, for example
in the $x_r$-\emph{chart}
$U_{x_r}$, we can choose coordinates $(y_1,\ldots,y_n)$ with
$y_1 = x_1,\ldots,y_r = x_r$, $y_{r+1} = x_{r+1}/x_r,\ldots,y_n = x_n/x_r$.
The following lemma is a simple but crucial exercise using the chain rule.

\begin{lemma} 
With coordinates chosen as above, we have:
\begin{align*}
\frac{\p}{\p x_j} &= \frac{\p}{\p y_j}, \qquad j = 1, \dots, r-1,\\
x_r \frac{\p}{\p x_r} &= y_r \frac{\p}{\p y_r} 
- \sum_{j=r+1}^n y_j \frac{\p}{\p y_j},\\
x_j \frac{\p}{\p x_j} &= y_j \frac{\p}{\p y_j}, \qquad j = r+1, \dots, n.
\end{align*}
\end{lemma}

The transformation formulas in Lemma 3.1 go back to \cite{giraud},
and are usually written a little differently. As arranged here, they
show that it is logarithmic differential operators which transform
naturally by blowing up.

\subsection{Derivative ideals}
Let $\cI \subset \cO_N$ denote a (coherent) ideal. Let $\cD(\cI) \subset
\cO_N$ denote the coherent sheaf of ideals generated by all first derivatives
of local sections of $\cI$ (so that $\cI \subset \cD(\cI)$); i.e., $\cD(\cI)$ is
the image of the natural morphism $\Der_N \times \cI \to \cO_N$.

If $f_1,\ldots,f_q$ are local generators of $\cI$ in a chart with
coordinates $(x_1,\ldots,x_n)$, then $f_i$, $\p f_i/\p x_j$, where 
$i = 1,\ldots,q$, $j= 1,\ldots,n$, are local generators of $\cD(\cI)$.

We also define \emph{higher-derivative ideals} inductively by
$$
\cD^{j+1}(\cI) := \cD(\cD^j(\cI)), \quad j= 1,\ldots.
$$

Now let $E$ denote a normal-crossings divisor on $N$. We define
$\cD_E(\cI) \subset \cO_N$ as the ideal generated by all local
sections of $\cI$ and first derivatives which preserve the ideal $\cI_E$
of $E$. (If $N \subset M$, where $M$ is smooth, and $E$ is a normal-crossings
divisor on $M$ which is transverse to $N$, then we will simply write
$\cD_E(\cI)$ instead of $\cD_{E|_N}(\cI)$.)

Consider a divisor $E = \sum_{i=1}^r \al_i H_i$, where each $\al_i \in \IN$.
Each point of $N$ has a neighbourhood $U$ with a coordinate system
$(x_1,\ldots,x_q,\ldots,x_n)$, $q \leq r$, such that $x_1,\ldots,x_q$
generate the ideals $\cI_H|_U$ for those $H = H_i$ which
intersect $U$; write $x_j = x_H$, where $H$ is the corresponding divisor.
If $f_i$ are local generators of $\cI$, then $\cD_E(\cI)$ is generated
locally by all
$$
f_i,\quad x_H\frac{\p f_i}{\p x_H},\quad  \mbox{and}\quad
\frac{\p f_i}{\p x_j},\quad
q+1 \leq j \leq n.
$$

We again define
$$\cD^{j+1}_E(\cI) := \cD_E(\cD^j_E(\cI)), \quad j= 1,\ldots.
$$

Now consider a marked ideal $\ucI = (M,N,E,\cI,d)$. For simplicity, we
will write $\ucI = (\cI,d)$ when the remaining entries are evident.
We define
\begin{align*}
\ucD(\ucI) &:= (M,N,E,\cD(\cI),d-1) = (\cD(\cI),d-1),\\
\ucD_E(\ucI) &:= (M,N,E,\cD_E(\cI),d-1) = (\cD_E(\cI),d-1).
\end{align*}
and
$$
\ucD_E^{j+1}(\ucI) := \ucD_E(\ucD_E^j(\ucI)) = (\cD^{j+1}_E(\cI),d-j-1),\quad
j= 1,\ldots.
$$

\begin{lemma}
Let $\ucI = (M,N,E,\cI,d)$ be a marked ideal. Then:
\begin{enumerate}
\item 
$\cD_E^k(\cD_E^l(\cI)) = \cD_E^{k+l}(\cI)$.
\item
$\cD_E^k(\cI\cdot\cJ) \subset \sum_{j=0}^k \cD_E^j(\cI)\cdot\cD_E^{k-j}(\cJ)$.
\item
$\cosupp \ucI \subset \cosupp \ucD_E^j(\ucI)$, $j=0,1,\ldots,$ with equality
if $E = \emptyset$.
\item
If $\vp: P \to N$ is a smooth morphism, then $\cD_{\vp^{-1}(E)}(\vp^*(\cI))
= \vp^*(\cD_E(\cI))$.
\end{enumerate}
\end{lemma}

\begin{proof}
These statements are immediate from the definitions. (It is enough to
check (4) on completions.)
\end{proof}

\begin{lemma}
Let $\ucI = (M,N,E,\cI,d)$ be a marked ideal. Let $\s: M' \to M$ denote
a test transformation for $\ucI$ 
(i.e., either an admissible or exceptional blowing-up,
or a projection $M \times \IA^1 \to M$), and let $\ucI' = (M',N',E',\cI',d')$
denote the transform of $\ucI$ by $\s$. (See Sect.\,2.) Then $\s$ is a
test transformation for $\ucD_E(\ucI)$, and $\cD_E(\cI)' \subset \cD_{E'}(\cI')$.
\end{lemma}

\begin{proof}
The first assertion follows from Lemma 3.2(3), and the second is an exercise
using the transformation formulas of Lemma 3.1.
\end{proof}

The following three corollaries are easy consequences of Lemma 3.3.

\begin{corollary}
With the hypotheses of Lemma 3.3,
if $j \leq d-1$, then $\cD^j_E(\cI)' \subset \cD^j_{E'}(\cI')$.
\end{corollary}

\begin{corollary}
If $j \leq d-1$, then any sequence of test transformations of $\ucI$ is
a sequence of test transformations of $\ucD^j_E(\ucI)$ and, after any sequence
of test transformations
\begin{equation*}
M = M_0 \stackrel{\s_1}{\longleftarrow} M_1 \longleftarrow \cdots
\stackrel{\s_{t}}{\longleftarrow} M_{t}\,,
\end{equation*}
we have an inclusion
$\cD^j_E(\cI)_t \subset \cD^j_{E_t}(\cI_t)$.
\end{corollary}

Of course, we can also define $\ucD^{j+1}(\ucI)$, $j=1,\ldots$, by iterating
$\ucD$, and we have analogues of Lemmas 3.2, 3.3 and Corollaries 3.4, 3.5
for the $\ucD^j(\ucI)$. We will use the following to prove Corollary 1.4. 

\begin{corollary} Suppose that $E= \emptyset$. Under the hypotheses of
Lemma 3.3, if $j \leq d-1$, then $\cD^j(\cI)' \subset \cD^j_{E'}(\cI')$.
Likewise after any sequence of test transformations, as in Corollary 3.5.
\end{corollary}

\subsection{Aside on sums and products of marked ideals} We define sums
and products, and give two simple lemmas.

Consider marked ideals
$\ucI = (M,N,E,\cI,d) = (\cI,d)$ and
$\ucJ = (M,N,E,\cJ,d) = (\cJ,d)$.
Define $\ucI\cdot\ucJ := (\cI\cdot\cJ,d+e)$.

\begin{lemma}
\begin{enumerate}
\item 
$\cosupp \ucI \cap \cosupp \ucJ \subset \cosupp \ucI\cdot\ucJ$
(whereas $\cosupp \ucI^k = \cosupp \ucI$.)
\item
Let $\s: M' \to M$ be a test morphism for both $\ucI$ and $\ucJ$. Then
$\s$ is a test morphism for $\ucI\cdot\ucJ$ and the transforms satisfy
$\ucI'\cdot\ucJ' = (\ucI\cdot\ucJ)'$.
\item
Multiplication of marked ideals is associative.
\end{enumerate}
\end{lemma}

Define $\ucI + \ucJ := (\ucI^{l/d} + \ucJ^{l/e}, l)$, where $l = \lcm(d,e)$.
Likewise, for any finite sum. Addition is not associative, but
$\ucI + \ucJ$ is equivalent to $(\ucI^e + \ucJ^d, de)$. (See Lemma 3.8(3).)

\begin{lemma}
\begin{enumerate} 
\item
$\cosupp (\ucI + \ucJ) = \cosupp \ucI \cap \cosupp \ucJ$.
\item
$\s: M' \to M$ is a test morphism for $\ucI + \ucJ$ if and only if
$\s$ is a test morphism for both $\ucI$ and $\ucJ$, and the transforms
by such a test morphism satisfy $\ucI' + \ucJ' = (\ucI + \ucJ)'$.
\item
Addition of marked ideals is associative up to equivalence.
\end{enumerate}
\end{lemma}

\subsection{Logarithmic derivatives determine equivalent ideals}
Consider a marked ideal $\ucI = (M,N,E,\cI,d)$. We define
$$
\ucC_E^k(\ucI) := \sum_{j=0}^k \ucD_E^j(\ucI), \quad k \leq d-1.
$$
Write $\ucC_E^k(\ucI) = \left(M,N,E, \cC_E^k(\ucI), d_{\ucC_E^k(\ucI)}\right)$.

\begin{lemma}
$\cosupp \ucI = \cosupp \ucC_E^k(\ucI).$
\end{lemma}

\begin{proof}
$\cosupp \ucI \subset \cosupp \ucC_E^k(\ucI)$, by Lemmas 3.2(3) and 3.8(1).
But $\cosupp \ucC_E^k(\ucI) \subset \cosupp \ucI$, since $\ucD_E^0(\ucI) = \ucI$.
\end{proof}

\begin{theorem}
Consider any sequence of test morphisms for $\ucI$,
\begin{equation}
M = M_0 \stackrel{\s_1}{\longleftarrow} M_1 \longleftarrow \cdots
\stackrel{\s_{t}}{\longleftarrow} M_{t}\,.
\end{equation}
Then (3.1) is a sequence of test morphisms for $\ucC_E^k(\ucI)$, and
for the transforms, we have
$$
\cosupp \ucC_{E_t}^k(\ucI_t) = \cosupp \ucC_E^k(\ucI)_t.
$$
\end{theorem}

\begin{proof}
By Corollary 3.5, any sequence of test morphisms of $\ucI$ is a sequence
of test morphisms of $\ucC_E^k(\ucI)$, and
$\cosupp \ucC_{E_t}^k(\ucI_t) \subset \cosupp \ucC_E^k(\ucI)_t$. Since
$\ucC_E^k(\ucI) = \ucI + \sum_{j=1}^k \ucD_E^j(\ucI)$, we have
\begin{align*}
\cosupp \ucC_E^k(\ucI)_t &\subset \cosupp \ucI_t\\
      &\subset \bigcap_{j=0}^k \cosupp \ucD_{E_t}^j(\ucI_t),\quad 
		 \mbox{by Lemma 3.2(3)},\\
      &= \cosupp \ucC_{E_t}^k(\ucI_t).
\end{align*}
\end{proof}

\begin{corollary}
$\ucI$ and $\ucC_E^k(\ucI)$ are equivalent.
\end{corollary}

\begin{proof}
They have the \emph{same} test sequences, by Lemma 3.9 and Theorem 3.10.
\end{proof}

\begin{lemma}
If $\vp: M' \to M$ is smooth, then
$\vp^*(\ucC_E^k(\ucI)) = \ucC_{\vp^{-1}(E)}^k(\vp^*(\ucI))$.
\end{lemma}

\begin{proof}
This follows from Lemma 3.2(4).
\end{proof}

\section{Maximal contact and coefficient ideals}

\begin{corollary}
Let $\ucI = (M,N,E,\cI,d)$ denote a marked ideal.
Let $z$ denote a section of $\cD_E^{d-1}(\cI)$ on $N|_U$, where $U$ is
an open subset of $M$. Suppose that $z$ has maximum order $1$ and that
$z$ is transverse to $E$ (equivalently, $z\cdot \cI_{E}|_N$ has only
normal crossings and $\cI_{E}|_N \not\subset (z)$, where $(z) \subset
\cO_N|_U$ is the ideal generated by $z$). Let $P = V(z)$. Then, after
any sequence of test transformations of $\ucI|_U$,
$$
\cosupp \ucI_t \subset P_t \subset N_t,
$$
where $P_t = V((z)_t)$ and $(z)_t$ denotes the transform of $(z)$
(cf. notation of Thm.\,3.10). Moveover, if $\cosupp \ucI \subsetneq P$,
then $(z)_t$ has maximum order $1$ and is transverse to $E_t$.
\end{corollary}

\begin{proof}
This follows from Theorem 3.10 and Corollary 3.11.
\end{proof}

In the situation of Corollary 4.1, we call the hypersurface $P$ of 
$N|_U$ a \emph{hypersurface of maximal contact} for $\ucI$. Note that
the existence of $z$ with maximum order $1$ means that $\ucI$ is of
maximal order on $U$. In the setting of Corollary 3.1, we define
the \emph{coefficient (marked) ideal}
\begin{equation}
\ucC_{E,P}(\ucI) := \left(U,P,E, \cC_E^{d-1}(\ucI)|_P, 
d_{\ucC_E^{d-1}(\ucI)}\right).
\end{equation}

\begin{corollary}
$\ucC_{E,P}(\ucI)$ is equivalent to $\ucI|_U$.
\end{corollary}

\begin{remark}
$\ucC_{E,P}$ commutes with smooth morphisms. In other words, if
$\vp: M' \to M$ is smooth, then $P' := \vp^{-1}(P) = V(\vp^*(z))$ is a
maximal contact hypersurface for $\vp^*(\ucI)$, and 
$\vp^*(\ucC_{E,P}(\ucI)) = \ucC_{E',P'}(\vp^*(\ucI))$ (by Lemma 3.12).
\end{remark}

\begin{exercises}
The following exercises give alternative ways to define coefficient
ideals. These approaches are used in \cite{BMjams}, \cite{BMinv}, 
and can simplify explicit calculations.
\smallskip

\noindent
(1)\ Under the hypotheses of Corollary 4.1, we can find a system of local
coordinates $(x_1,\ldots,x_n)$ for $N$ such that $x_n = z$ and the
components of $E$ are given by $x_i = 0$, $i = 1,\dots,r<n$ (in this
chart). Let $\cD_z(\cI)$ denote the ideal generated by $f$, $\p f/\p z$,
for all $f \in \cI$ (so that $\cI \subset \cD_z(\cI) \subset \cD_E(\cI)$),
and set $\ucC_z^{d-1}(\ucI) := \sum_{j=0}^{d-1}\ucD_z^j(\ucI)$. Then
$\ucC_{z,P}(\ucI) := \left(U,P,E, \cC_z^{d-1}(\ucI)|_P, 
d_{\ucC_z^{d-1}(\ucI)}\right)$ is equivalent to $\ucI$ (in the chart).
\smallskip

\noindent
(2)\ Suppose that $\cI \subset \cO_M$ is a principal ideal. Let $g$ be a
local generator of $\cI$ in a neighbourhood $U$ of $a \in M$ with coordinates
$(x_1,\ldots,x_m)$. Let $d = \ord_a\cI$.
(By making a linear coordinate change), we can assume that
$\p^d g/\p x_m^d$ is nonvanishing. Let $z := \p^{d-1} g/\p x_n^{d-1}$.
Then $N := V(z)$ is a hypersurface of maximal contact 
for $\ucI|_U := (U,U,\emptyset, \cI|_U, d)$, and $\ucI|_U$ is equivalent
to $\ucC_{x_n,N}(\ucI) := \left(U,N,\emptyset, \cC_{x_n}^{d-1}(\cI)|_N,
d_{\ucC_{x_n,N}(\ucI)}\right)$, where $\cC_{x_n}^{d-1}(\cI) :=
\sum_{k=0}^{d-1} \left((\p^k g/\p x_n^k), d-k\right)$.
\end{exercises}

\section{Proof of resolution of singularities of a marked ideal}
In this section, we prove Theorem 1.3.
Let $\ucI = (M,N,E,\cI,d)$ denote a marked ideal, where $d>0$. The proof is by
induction on $\dim \ucI := \dim N$.

First suppose $\dim N = 0$. If $\cI=0$, then $\cosupp \ucI = N$ and we
can blow up with centre $N$ to resolve the singularities of $\ucI$.
If $\cI \neq 0$, then $\cosupp \ucI = \emptyset$ (since $d>0$), so that
$\ucI$ is already resolved.

The inductive step breaks up into the two independent steps I and II
formulated in Section 1.
\medskip

\noindent
{\bf Step I.} Let $\ucI = (M,N,E,\cI,d)$ be a marked ideal of maximal order.
\smallskip

\noindent
{\bf Case A. $E = \emptyset$.} Let $a \in \cosupp \ucI$. By Corollary 4.1,
there is a hypersurface of maximal contact $P$ for $\ucI$ at $a$; $P \subset
N\cap U$, where $U$ is a neighbourhood of $a$ in $M$. Define
\begin{equation}
\ucC_U(\ucI) := \ucC_{\emptyset,P}(\ucI) 
= (U,P,\emptyset, \ucC_{\emptyset}^{d-1}(\ucI)|_P, d_{\ucC}),
\end{equation}
where $d_{\ucC} = d_{{\ucC}_{\emptyset}^{d-1}(\ucI)}$. Then $\ucC_U(\ucI)$ is
equivalent to $\ucI|_U = (U,N\cap U, \emptyset, \cI|_{N\cap U}, d)$,
according to Corollary 4.2. Therefore, a resolution of singularities of
$\ucC_U(\ucI)$ (which exists by induction) is a resolution of 
singularities of $\ucI|_U$.

The marked ideals $\ucC_U(\ucI)$ and $\ucC_V(\ucI)$ defined in two overlapping
charts $U$ and $V$ are equivalent in $U\cap V$ (since both are equivalent
to $\ucI_{U\cap V}$; therefore, by functoriality in dimension $n-1$, their
resolution sequences  are the same over $U\cap V$ (not counting
blowings-up that restrict to isomorphisms over $U\cap V$). 

\begin{claim}
It follows that the resolution sequences for the $\ucC_U(\ucI)$ (defined
locally) patch together to give functorial resolution of singularities
of a marked ideal $\ucI$.
\end{claim}

The problem here is that, to glue together the local centres of blowing up,
we have to know which nontrivial centres in the various charts $U$ should be
taken first, etc. We leave this claim to Section 7 below, where we give
two different proofs.
\smallskip

\noindent
{\bf Case B. General maximal order case.} Let $\ucI_{\emptyset}$ denote
the marked ideal $(M,N, \emptyset, \cI,d)$. Let $a \in \cosupp \ucI
= \cosupp \ucI_{\emptyset}$. By Corollary 4.1, there exists an (irreducible)
hypersurface of maximal contact $P$ for $\ucI_{\emptyset}$, defined in a
neighbourhood $U$ of $a$. We introduce $\ucC := \ucC_U(\ucI_{\emptyset})
= (U,P,\emptyset, \cC, d_{\ucC})$, as in Case A above.

If $x \in N$, set 
$$
s(x) := \#\{H \in E: x \in H\}
$$
(where ``$H \in E$'' means ``$H$ is a component of $E$'', and $\#$
denotes the cardinality of a finite set). After any sequence of test 
transformations of $\ucI$, we will continue to write $E$ (as opposed to $E_t$)
for the divisor whose components are the strict transforms of those of $E$,
with the same ordering as in $E$, and we will continue to write $s(x)$ for
$\#\{H \in E: x \in H\}$.

Let $s := \max \{s(x): x \in \cosupp \ucI\}$,
and let $\Sub(E,s)$ denote the set of $s$-element subsets of $E$
(i.e., of the set of components of $E$). Let $\ucE$ denote the marked ideal
$\ucE := \left(U,P,\emptyset, \cE, d_{\ucC}\right)$, where
$$
\cE := \prod_{\La \in \Sub(E,s)}\sum_{H \in \La} \cI_H^{d_\cC}\cdot \cO_P.
$$
Set $\ucJ := \ucC + \ucE$. 

Clearly, $\cosupp \ucE = \{x \in P: s(x) \geq s\}$, so that
$$
\cosupp \ucJ = \cosupp \ucI|_U \bigcap \{x \in P: s(x) = s\},
$$
(and likewise after any sequence of test transformations of $\ucJ$).
Thus any sequence of test transformations of $\ucJ$ is a sequence
of test transformations of $\ucI|_U$ (and the centre $C$ of each admissible 
blowing-up lies inside the intersection of the components of $E$
that pass through a point of $C$). Therefore, the equivalence class
of $\ucJ$ depends only on that of $\ucI|_U$, and blowings-up which are 
admissible for $\ucJ$ are also admissible for $\ucI|_U$.

By induction on dimension, there is a resolution of singularities of $\ucJ$
We can define $\ucJ$ for each element of a finite covering $\{U\}$ of $M$
and, as in Claim 4.1, the resolution sequences for the $\ucJ$ patch together
to define a sequence of admissible blowings-up of $\ucI$. (See Section 7.)

If this sequence of blowings-up resolves $\ucI$, then we have finished
Case B. Otherwise, 
$\cosupp \ucI_j$ and $\{x: s(x)=s\}$ become disjoint
for some index $j$ of the blowings-up sequence, say for $j=q_1$ 
(i.e., $s(x) < s$, for all $x \in \cosupp \ucI_j$).

We can now repeat the process above, using $\ucC_{q_1} = 
\ucC_U(\ucI_{\emptyset})_{q_1}$, $\ucC_{q_2} = \ucC_U(\ucI_{\emptyset})_{q_2}$,
etc., in place of $\ucC$, and using the new value of $\max s(x)$ in place of $s$.
After finitely many steps, $\ucI$ is resolved, so we have completed Case B
and therefore Step I.
\medskip

\noindent
{\bf Step II. General case.} Let $\ucI = (M,N,E,\cI,d)$ be an arbitrary
marked ideal.

First suppose that $\cI = 0$ (so that $\cosupp \ucI = N$). Then we can
blow up with centre $N$ to resolve singularities.

Now suppose that $\cI \neq 0$. Write
\begin{equation}
\cI = \cM(\ucI)\cdot \cN(\ucI),
\end{equation}
where $\cM(\ucI)$ is a product of prime ideals defining irreducible components
of the elements of $E$, and $\cN(\ucI)$ is divisible by no such prime ideal.
We consider two cases.
\smallskip

\noindent
{\bf Case A. Monomial case $\cI = \cM(\ucI)$.} For any sequence of
admissible blowings-up $\s_j$ of $\ucI$, let us order the collection of subsets
of each $E_j$ as follows: Suppose that $E = \{H_1,\ldots, H_q\}$. Write
$E_0 := E = \{H_1,\ldots, H_q\} =: \{H_1^0,\ldots, H_q^0\}$. For each
$j=1,\ldots$, write $E_j = \{H_1^j,\ldots, H_{q+j}^j\}$, where $H_i^j$
denotes the strict transform of $H_i^{j-1}$ by $\s_j$, if $i < q+j$, and
$H_{q+j}^j = \s_j^{-1}(C_{j-1})$. Now, we associate to each subset $I$ of
$E_j$ the finite sequence $(\de_1,\ldots,\de_{q+j})$, where
$\de_i = 0$ if $H_i^j \not\in I$, and $\de_i = 1$ if $H_i^j \in I$.
Then we order the subsets $I$ of $E_j$, for \emph{all} $j$, by the 
lexicographic ordering of such sequences.

Let $a \in \cosupp \ucI$. In a neighbourhood of $a$, we can write
$\cI = \cI_{H_{i_1}}^{\al_1}\cdots \cI_{H_{i_r}}^{\al_r}$
(where $a \in H_{i_1} \cap \cdots \cap H_{i_r}$ and where we write $\cI_{H_{i_k}}$
instead of $\cI_{H_{i_k}}\cdot\cO_N$ to simplify the notation);
in particular,
$\al_1 + \cdots + \al_r > d$. (In Theorem 6.1 below, we will see
that $\mu_a(\ucI) = (\al_1 + \cdots + \al_r)/d$ is a local invariant 
of the equivalence
class of $\ucI$.) Then (in the neighbourhood above) 
$\cosupp \ucI = \cup Z_I$, where each $Z_I := N \cap \bigcap_{H \in I}H$,
and $I$ runs over the \emph{smallest} subsets of $\{H_{i_1},\ldots,H_{i_r}\}$
such that $\sum_{l \in I} \al_l \geq d$;
in other words, $I$ runs over the subsets of $\{H_{i_1},\ldots,H_{i_r}\}$
such that
$$
0 \leq \sum_{l \in I} \al_l -d < \al_k, \quad \mbox{for all } k \in I.
$$
(We have simplified the notation by identifying subsets of
$\{H_{i_1},\ldots,H_{i_r}\}$ with subsets of $\{1,\ldots,r\}$.)

Let $J(a)$ denote the maximum of these subsets $I$ (with respect to the
order above). Clearly, $J(x)$ is Zariski upper-semicontinuous on 
$\cosupp \ucI$, and the maximum locus of $J(x)$ consists of at most one
irreducible component of $\cosupp \ucI$ through each point. 

Consider the blowing-up $\s$ with centre $C$ given by the maximum locus
of $J(\cdot)$. If $a \in C$, then (in a neighbourhood as above),
$C = N\cap \bigcap_{l \in J(a)}H_{i_l}$. We can choose local coordinates
$(x_1,\ldots,x_n)$ for $N$ at $a$ such that, for each $k \in J(a)$,
$x_k$ is a local generator $x_{H_{i_k}}$ of the ideal $\cI_{H_{i_k}}|_N$.
Then, in the $x_k$-chart of the blowing-up $\s$, the transform of $\cI$
is given by
$$
\cI' = \cI_{H_{i_1}'}^{\al_1}\cdots \cI_{\s^{-1}(C)}^{\sum_{J(a)}\al_l-d}
\cdots \cI_{H_{i_r}'}^{\al_r},
$$
(with $\cI_{\s^{-1}(C)}^{\sum_{J(a)}\al_l-d}$ in the $k$'th place). Since
$\sum_l \al_l -d < \al_k$, $\mu_{a'}(\ucI') = \mu_a(\ucI) - p/d$, where
$p$ is a positive integer, over all points $a$ of some component of
$\cosupp \ucI$. We therefore resolve singularities after a finite number
of steps.
\smallskip

\noindent
{\bf Case B. General case.} Set
$$
\ord\,\cN(\ucI) := \max_{x \in \cosupp \ucI} \ord_x \cN(\ucI),
$$
and define marked ideals 
\begin{align*}
\ucN(\ucI) &:= (\cN(\ucI), \ord\,\cN(\ucI)),\\ 
\ucM(\ucI) &:= (\cM(\ucI), d-\ord\,\cN(\ucI)).
\end{align*}
We define the \emph{companion ideal} of $\ucI = (\cI,d)$ as
$$
\ucG(\ucI) :=
\begin{cases}
\ucN(\ucI) + \ucM(\ucI), &\text{if $\ord\,\cN(\ucI) < d$;}\\
\ucN(\ucI),              &\text{if $\ord\,\cN(\ucI) \geq d$.}
\end{cases}
$$
Then $\ucG(\ucI)$ is of maximal order and, for any sequence of test
transformations of $\ucG(\ucI)$,
\begin{align*}
\cosupp \ucG(\ucI)_j &= \cosupp \ucN(\ucI)_j \bigcap \cosupp \ucM(\ucI)_j\\
		     &\subset \cosupp \ucN(\ucI)_j \bigcap \cosupp \ucI_j.
\end{align*}
By Lemma 3.7, $\ucI_j = \ucM(\ucI)_j\cdot\ucN(\ucI)_j$. 

\begin{remark}
For any sequence
of \emph{admissible} blowings-up of $\ucG(\ucI)$, it is clear that
$\ucN(\ucI)_j = \ucN(\ucI_j)$, $\ucM(\ucI)_j = \ucM(\ucI_j)$, and
$\cosupp \ucG(\ucI)_j = \cosupp \ucN(\ucI)_j \cap \cosupp \ucI_j$.
\end{remark}

In Theorems 6.1 and 6.2 below, we will prove that, if $a \in \cosupp \ucI$,
then 
\begin{equation}
\mu_a(\ucI) := \frac{\ord_a\cI}{d} \quad \mbox{and} \quad 
\mu_{H,a}(\ucI) := \frac{\ord_{H,a}\cI}{d},\,\, H \in E,
\end{equation}
depend only on the equivalence class of $\ucI$ and $\dim N$. 
($\ord_{H,a}\cI$ denotes the \emph{order} of $\cI \subset \cO_N$
\emph{along} $H|_N$ at $a$;\, i.e., the largest $\mu \in \IN$ such that
$\cI_a \subset \cI_{H|_N,a}^\mu$.) The following
is a corollary of these results.

\begin{corollary}
The equivalence class of $\ucG(\ucI)$ depends only on the equivalence
class of $\ucI$ and $\dim N$.
\end{corollary}

\begin{proof}
The divisor (with rational coefficients) $\cM(\ucI)^{1/d}$ depends only on
the equivalence class of $\ucI$, by Theorems 6.1 and 6.2. Over any sequence
of test transformations of $\ucG(\ucI)$,
\begin{alignat*}{2}
\cosupp \ucG(\ucI)_j &= \left\{x \in \cosupp \ucI_j:\right. 
&&\left.\ord_x\cN(\ucI)_j \geq \ord\,\cN(\ucI),\right.\\
&&&\left.\ord_x\cM(\ucI)_j \geq d - \ord\,\cN(\ucI)\right\}\\
&= \left\{x \in \cosupp \ucI_j:\right. 
&&\left.\mu_x(\ucI_j) - \frac{1}{d} \ord_x\cM(\ucI)_j \geq 
\frac{1}{d} \ord\,\cN(\ucI),\right.\\
&&&\left.\frac{1}{d} \ord_x\cM(\ucI)_j \geq 1 - \frac{1}{d} \ord\,\cN(\ucI)\right\}.
\end{alignat*}
All terms in the latter condition depend only on the equivalence class of 
$\ucI$ (and $\dim N$).
\end{proof}

Now, we can resolve the singularities of the marked ideal of maximal order
$\ucG(\ucI)$, using Step I. The blowings-up involved are admissible for $\ucI$.
Such a resolution leads, after a finite number of steps (say $r_1$ steps) to
a marked ideal $\ucI_{r_1}$ such that either $\ord\,\cN(\ucI_{r_1}) <
\ord\,\cN(\ucI)$, or $\cosupp \ucI_{r_1} = \emptyset$. (See Remark 5.2.)
In the latter case, we have resolved the singularities of $\ucI$. In the
former case, we can repeat the process for $\ucI_{r_1} = (\cI_{r_1},d),
\ldots$. We get marked ideals $\ucI_{r_1},\ldots,\ucI_{r_k},\ldots$, such that
$$
\ord\,\cN(\ucI) > \ord\,\cN(\ucI_{r_1}) > \cdots > \ord\,\cN(\ucI_{r_k}) > \cdots .
$$
The process terminates after a finite number of steps, when either we have 
$\ord\,\cN(\ucI_{r_k}) = 0$ (i.e., we have reduced to the monomial case
$\ucI_{r_k} = \ucM(\ucI_{r_k})$), or we have $\cosupp \ucI_{r_k} = \emptyset$
(i.e., we have resolved singularities).

This will complete the proof of functorial resolution of singularities
of a marked ideal (Theorem 1.3) once we establish Claim 5.1 and prove Theorems
6.1, 6.2. \hfill $\Box$

\section{Invariants of a marked ideal}
If $\ucI = (M,N,E,\cI,d)$ is a marked ideal and $a \in \cosupp \ucI$,
we define $\mu_a(\ucI) := \ord_a\cI/d$ and 
$\mu_{H,a}(\ucI) := \ord_{H,a}\cI/d$, $H \in E$, as in (5.3).
The following theorem is due to Hironaka \cite{Hid}.

\begin{theorem} Let $\ucI = (M,N,E,\cI,d)$ and $\ucJ = (M,P,E,\cJ,e)$
be marked ideals. Suppose that $\ucI$ and $\ucJ$ are weak-equivalent
(Defns.\,2.4). Let $a \in \cosupp \ucI = \cosupp \ucJ$. If $\dim N
= \dim P$, then $\mu_a(\ucI) = \mu_a(\ucJ)$. If $\dim N > \dim P$, then 
$\mu_a(\ucI) = 1$.
\end{theorem}

\begin{proof}
We have $\mu_a(\ucI) = \infty$ if and only if $\cI = 0$ at $a$, i.e.,
$\cosupp \ucI$ is smooth of dimension $= \dim N$ at $a$.

Suppose that $\mu_a(\ucI) < \infty$. Set $\ucI_0 := \ucI$, 
$\ucI_0 = (M_0, N_0, E_0, \cI_0, d)$. Let $\s_1: M_1 := M_0 \times \IA^1
\to M_0$ denote the projection, and let $\ucI_1 = (M_1, N_1, E_1, \cI_1, d)$
denote the transform of $\ucI_0$ by $\s_1$ (Defns.\,2.4). Let $a_1 :=
(a,0)$ and let $\G_1 := \{a\}\times \IA^1 \subset \cosupp \ucI_1 =
\s_1^{-1}(\cosupp \ucI_0)$. Consider a sequence of admissible blowings-up
\begin{equation*}
M_1 \stackrel{\s_2}{\longleftarrow} M_2 \longleftarrow \cdots
\longleftarrow M_j \stackrel{\s_{j+1}}{\longleftarrow} M_{j+1}
\longleftarrow \cdots
\end{equation*}
and the corresponding transforms $\ucI_j = (M_j, N_j, E_j, \cI_j, d)$
of $\ucI_1$, where, for each $j\geq 1$, $\s_{j+1}$ denotes the blowing-up
with center $a_j$, and the $a_j$ are defined inductively as follows:
For each $j\geq 1$, $a_{j+1} := \G_{j+1} \cap D_{j+1}$, where
$D_{j+1}$ is the exceptional divisor $\s_{j+1}^{-1}(a_j)$ and $\G_{j+1}$
is the strict transform of $\G_j$.

Clearly, each $\G_{j+1} \subset \cosupp \ucI_{j+1}$, since $\G_{j+1}\setminus
\{a_{j+1}\} \subset \cosupp \ucI_{j+1}$.

We introduce a subset $S$ of $\IN \times \IN$ depending only on $\dim N$
and the weak-equivalence class of $\ucI$, as follows: Let $r = \dim M -\dim N$.
Note that, for each $j \geq 1$, $D_{j+1} \cap N_{j+1}$ is smooth, and has
codimension $r$ in $D_{j+1}$ and codimension $1$ in $N_{j+1}$. We say that
$(j,0) \in S$, $j\geq 1$, if, after $j$ blowings-up $\s_2,\ldots,\s_{j+1}$
as above, $D_{j+1} \cap \cosupp \ucI_{j+1}$ contains (i.e., is) a smooth
subvariety of codimension $r$ in $D_{j+1}$ at $a_{j+1}$. 

In this case, $D_{j+1} \cap \cosupp \ucI_{j+1} = D_{j+1} \cap N_{j+1}$ at
$a_{j+1}$, and the condition $D_{j+1} \cap N_{j+1} \subset \cosupp \ucI_{j+1}$
means that $\ord_{D_{j+1},a_{j+1}}\ucI_{j+1} \geq d$. The latter inequality is
equivalent to $j(\mu_a(\ucI) - 1)d \geq d$, i.e., $j(\mu_a(\ucI) - 1) \geq 1$.
(This can be proved by induction on $j$, by a simple formal Taylor series
calculation: For each $j$, $\cI_{j+1,a_{j+1}} =
\cI_{D_{j+1},a_{j+1}}^{-jd}\cdot(\cI_a\cdot\cO_{N_{j+1},a_{j+1}})$. 
But $\ord_{a_{j+1}}\ucI_{j+1}
= \ord_{a_j}\ucI_j = \cdots = \ord_a\ucI$; therefore, 
$\cI_a\cdot\cO_{N_{j+1},a_{j+1}}$ is divisible by 
$\cI_{D_{j+1},a_{j+1}}^{j\ord_a\cI}$ and not by any higher power. Since
$\ord_a\cI = \mu_a(\ucI)d$, we have $\cI_{j+1,a_{j+1}} = 
\cI_{D_{j+1},a_{j+1}}^{j(\mu_a(\ucI)-1)d}\cdot\tilde{\cI}_{a_{j+1}}$,
where $\tilde{\cI}_{a_{j+1}}$ is not divisible by $\cI_{D_{j+1},a_{j+1}}$.)

In particular, since $\mu_a(\ucI)\geq 1$, $(j,0) \not\in S$ for all $j\geq 1$
if and only if $\mu_a(\ucI)=1$.

It follows that if $\ucI$ is weak-equivalent to $\ucJ$, where $\dim P < \dim N$,
then $\mu_a(\ucI) = 1$ (since $D_{j+1} \cap \cosupp \ucI_{j+1} = 
D_{j+1} \cap \cosupp \ucJ_{j+1} \subset D_{j+1} \cap P_{j+1}$ cannot contain a
smooth subvariety of codimension $r$ in $D_{j+1}$.)

Now suppose that $(j,0) \in S$, for some $j \geq 1$. Then we can consider
the blowing-up $\s_{j+2}: M_{j+2} \to M_{j+1}$ with centre $C_{j+1} =
D_{j+1} \cap N_{j+1}$. Note that $\s_{j+2}|_{N_{j+2}}: N_{j+2} \to N_{j+1}$
is the identity! Set $G_{j+2} := \mbox{strict transform of } G_{j+1}$ and
$a_{j+2} := G_{j+2} \cap D_{j+2}$. We say that $(j,1) \in S$ if 
$D_{j+2} \cap \cosupp \ucI_{j+2}$ contains a smooth subvariety of codimension
$r$ in $D_{j+2}$ at $a_{j+2}$.

If so, then again $D_{j+2} \cap \cosupp \ucI_{j+2} = D_{j+2} \cap N_{j+2}$.
Since $\cI_{j+1} = \cI_{D_{j+1}}^{j(\mu_a(\ucI)-1)d}\cdot\tilde{\cI}_{j+1}$
at $a_{j+1}$, where $\tilde{\cI}_{j+1}$ is not divisible by $\cI_{D_{j+1}}$, we
see that $\cI_{j+2} = \cI_{D_{j+2}}^{j(\mu_a(\ucI)-1)d-d}\cdot\tilde{\cI}_{j+2}$
at $a_{j+2}$, where $\tilde{\cI}_{j+2} = \tilde{\cI}_{j+1}$ is not divisible by
$\cI_{D_{j+2}}$. Therefore, $(j,1) \in S$ if and only if $j(\mu_a(\ucI)-1)d-d
\geq d$.

We continue inductively: If $i\geq 1$ and $(j,i-1) \in S$, let
$\s_{j+i+1}: M_{j+i+1} \to M_{j+i}$ denote the blowing-up with centre
$C_{j+i} = D_{j+i} \cap N_{j+i}$. Set $G_{j+i+1} := \mbox{ strict transform of }
G_{j+i}$ and $a_{j+i+1} := G_{j+i+1} \cap D_{j+i+1}$. We say that $(j,i) \in S$ if
$D_{j+i+1} \cap \cosupp \ucI_{j+i+1}$ is smooth and of codimension $r$ in 
$D_{j+i+1}$; i.e., $(j,i) \in S$ if and only if $j(\mu_a(\ucI)-1) - i \geq 1$.

According to its definition, $S$ depends only on $\dim N$ and the weak-equivalence
class of $\ucI$ (and the point $a$). We have shown that $S = \emptyset$
if and only if
$\mu_a(\ucI) = 1$, and, if $S \neq \emptyset$, then 
$$
S = \left\{(j,i)\in \IN \times \IN: j(\mu_a(\ucI)-1) - i \geq 1\right\};
$$
therefore, $\mu_a(\ucI)$ is uniquely determined by $S$. ($\mu_a(\ucI) =
1 + \sup_{(j,i)\in S}(i+1)/j$.)
\end{proof}

\begin{theorem}
Let $\ucI = (M,N,E,\cI,d)$ be a marked ideal and let $a \in \cosupp \ucI$.
Consider $H \in E$. Then $\mu_{H,a}(\ucI)$ depends only
on the equivalence class of $\ucI$ and the dimension of $N$ at $a$. (See
Defns.\,2.5.)
\end{theorem}

\begin{proof}
Suppose that $H \ni a$. Set Set $\ucI_0 := \ucI$,
$\ucI_0 = (M_0, N_0, E_0, \cI_0, d)$. We can choose local coordinates
$(x_1,\ldots,x_n)$ for $N$ at $a$, such that $H|_{N_0} = \{x_1=0\}$. Let $\s_1:
M_1 = M_0 \times \IA^1 \to M_0$ denote the projection, let $a_1 = (a,0)$
and let $\ucI_1 = (M_1, N_1, E_1, \cI_1, d)$ denote the transform of 
$\ucI_0$ by $\s_1$. (Thus $(x_1,\ldots,x_n,t)$ are local coordinates for
$N_1$ at $a_1$.) Let $H_0^1 := M_0 \times \{0\}$ (so that $H_0^1|_{N_1}
= \{t=0\}$) and let $H_1^1 := \s^{-1}(H)$ ((so that $H_1^1|_{N_1} = \{x_1=0\}$.)
Set $\G_1 := \{a\} \times \IA^1$.

We follow $\s_1$ by a sequence of exceptional blowings-up
\begin{equation*}
M_1 \stackrel{\s_2}{\longleftarrow} M_2 \longleftarrow \cdots
\longleftarrow M_j \stackrel{\s_{j+1}}{\longleftarrow} M_{j+1}
\longleftarrow \cdots,
\end{equation*}
where (for all $j\geq 1$) we inductively define $\s_{j+1}$ as the blowing-up
with centre $C_j := H_0^j \cap H_1^j$, and set $H_0^{j+1} :=
\s_{j+1}^{-1}(C_j)$, $H_1^{j+1} := \mbox{strict transform of }
H_1^j$ by $\s_{j+1}$, $\G_{j+1} := \mbox{ strict transform of } \G_j$,
$a_{j+1} := \G_{j+1} \cap H_0^{j+1} = C_{j+1} \cap \s_{j+1}^{-1}(a_j)$. 
Then, for each $j\geq 1$, $N_{j+1}$ has a chart
with coordinates $(x,t) = (x_1,\ldots,x_n,t)$ at $a_{j+1}$ in which
$(\s_2\circ \cdots \circ \s_{j+1})(x,t) = (t^jx_1,x_2, \ldots, x_n,t)$.

Write $\cI = \cI_H^{\mu_{H,a}(\ucI)d}\cdot\cJ$ at $a$, where
$\cJ$ is not divisible by $\cI_H$. Set $\tau_{j+1}:= \s_1\circ\s_2\circ
\cdots\circ\s_{j+1}$. Then $\tau_{j+1}^*(\cJ)$ is not divisible by 
$\cI_{H_1^{j+1}}$ at $a_{j+1}$, and 
$$
\tau_{j+1}^*(\cI) = \cI_{H_0^{j+1}}^{j\mu_{H,a}(\ucI)d}\cdot 
\cI_{H_1^{j+1}}^{\mu_{H,a}(\ucI)d}\cdot \tau_{j+1}^*(\cJ).
$$
(Consider the formal expansion of the pullback of $f(x)=x_1^{\mu_{H,a}(\ucI)d}g(x)
\in \cI$.) 

Therefore, there exists $i \geq 1$ such that $\ord_{a_{i+1}} \tau_{i+1}^*(\cJ)
= \ord_{a_{j+1}} \tau_{j+1}^*(\cJ)$, for all $j \geq i$. (We can take $i$ to be
the least order of a monomial not divisible by $x_1$ in the Taylor expansions
of a set of generators of $\cJ$ at $a$.) Clearly, if $j\geq i$, then
$$
\mu_{H,a}(\ucI) = \mu_{a_{j+1}}(\ucI_{j+1}) - \mu_{a_j}(\ucI_j),
$$
so the result follows from Theorem 6.1.
\end{proof}

\section{Passage from local to global}
In this section, we give two different proofs of Claim 5.1 in Step I.A in
Section 5 (and of the analogous claim in Step I.B). The first proof is the
method described by Kollar \cite{Ko}, which we present briefly here. (See
\cite[Prop. 98]{Ko} for an axiomatization of the procedure.) The second proof
is the method of \cite{BMinv, BMda1}, which involves proving a stronger result 
than Theorem 1.3: we give an algorithm for functorial resolution of singularities,
where each successive centre of blowing up is the maximum locus of a 
\emph{desingularization invariant}. We use the notation of Section 4.

The invariant is essentially already present in the proof of Theorem 1.3 in
Section 5; we just make it explicit below. The invariant is defined using a
sequence of pairs --- it begins with $(\ord\,\cN(\ucI), s)$ and the following
pairs are determined in decreasing dimension, by induction. (The first pair
is the same as that denoted $(d,s)$ in our first version of resolution of
singularities \cite[\S4]{BMihes}.)

\subsection{First proof} For simplicity, assume that $N$ has constant dimension $n$.
Choose a finite open covering $M = \cup U_i$ such that, for each $i$, we have
a coefficient ideal $\ucC_{U_i}(\ucI)$ as in Step I.A. Set $M' := \coprod_i U_i$
(the disjoint union) and let $\g: M' \to M$ denote the natural morphism. ($\g$
is smooth.) The fibre-product $M'' := M' \times_M M'$ can be identified with the
disjoint union $\coprod_{i \leq j} U_i \cap U_j$, so that, if $\tau_1,\ \tau_2:
M'' \to M'$ are the two projections, then $\tau_1|_{U_i \cap U_j}$ and 
$\tau_2|_{U_i \cap U_j}$ are the inclusions $U_i \cap U_j \hookrightarrow U_i$
and $U_i \cap U_j \hookrightarrow U_j$ (respectively).

Set $\ucI' := \g^*(\ucI)$ and $\ucI'' := \tau_1^*(\ucI') = \tau_2^*(\ucI')$.
The coefficient ideals $\ucC_{U_i}(\ucI)$ induce a marked ideal $\ucC'$ of
dimension $n-1$ on $M'$. Set $\ucC_1'' := \tau_1^*(\ucC')$, $\ucC_2''
:= \tau_2^*(\ucC')$. Then $\ucC'$ is equivalent to $\ucI'$, and $\ucC_1''$ is 
equivalent to $\ucC_2''$.

By the inductive assumption of functorial resolution of singularities in
dimension $n-1$, there are resolution sequences associated to $\ucC'$,
$\ucC_1''$ and $\ucC_2''$ such that: (1) the resolution sequences for
$\ucC_1''$ and $\ucC_2''$ are the same; (2) the morphism $\tau_1$ (or $\tau_2$)
lifts throughout the resolution sequences for $\ucC_1''$ and $\ucC'$
(or $\ucC_2''$ and $\ucC'$).

Since $\ucC'$ is equivalent to $\ucI'$, the resolution sequence for 
$\ucC'$ is a resolution sequence for $\ucI'$. Consider the centre $C_0' \subset
M_0' := M'$ of the first blowing-up in this sequence. For each $i$, let $C_{0i}'$
denote the restriction of $C_0'$ to $U_i$. Then $C_{0i}'|{U_i \cap U_j} =
C_{0j}'|{U_i \cap U_j}$, for all $i \leq j$.

Therefore, the $C_{0i}'$ glue together to define a smooth closed subvariety
$C_0$ of $M_0 := M$. We thus obtain the centre of the first global blowing-up,
and can continue
in the same way to get the entire functorial resolution sequence for $\ucI$.

The same argument works for the analogous claim in Step I.B.

\subsection{Second proof} Given a marked ideal $\ucI$, we will construct
a resolution sequence and an invariant (with values in an ordered set)
defined at the points of the 
$\cosupp \ucI_j$ (where the $\ucI_j$ are the successive transforms of $\ucI$),
such that the invariant is upper-semicontinuous, we blow up its
maximum locus at each step, and the invariant decreases with (finitely many)
blowings-up.
The following theorem makes this precise. Let $\IQ_{>0}$ (or $\IQ_{\geq 0}$)
denote the set of positive (or nonnegative) rational numbers.

\begin{theorem} Let $\ucI = (M,N,E,\cI,d)$ denote a marked ideal, where $d>0$.
Then $\ucI$ admits a resolution of singularities 
\begin{equation*}
M = M_0 \stackrel{\s_1}{\longleftarrow} M_1 \longleftarrow \cdots
\stackrel{\s_{t}}{\longleftarrow} M_{t}\,,
\end{equation*}
such that, if $\ucI_j$ denotes the $j$'th transform of $\ucI$, then the 
following properties are satisfied.
\begin{enumerate}
\item
There are Zariski upper-semicontinuous functions $\inv_{\ucI}$, $\mu_{\ucI}$
and $J_{\ucI}$ defined on $\cosupp \ucI_j$, for all $j$, where 
$\inv_{\ucI}$ takes values in the set of sequences consisting of finitely many
pairs in $\IQ_{>0}\times \IN$ followed by $0$ or $\infty$ (ordered 
lexicographically), $\mu_{\ucI}$ takes values in $\IQ_{\geq 0} \cup \{\infty\}$,
and $J_{\ucI}$ takes values in the set of subsets of $E_j$, for all $j$ (ordered
as in Step II.A in Section 5).
\item
Each centre of blowing up $C_j \subset M_j$ is given by the maximum locus
of $(\inv_{\ucI}, J_{\ucI})$ on $\cosupp \ucI_j$ (where pairs are ordered
lexicographically).
\item
Let $a \in \cosupp \ucI_{j+1}$ and $b = \s_{j+1}(a)$. If $b \not\in C_j$, then
$$
\inv_{\ucI}(a)=\inv_{\ucI}(b),\quad \mu_{\ucI}(a)=\mu_{\ucI}(b),\quad
J_{\ucI}(a)=J_{\ucI}(b).
$$
If $b \in C_j$, then
$$
\left(\inv_{\ucI}(a),\mu_{\ucI}(a)\right) 
< \left(\inv_{\ucI}(b),\mu_{\ucI}(b)\right).
$$
\item
Suppose that $\vp: M' \to M$ is a smooth morphism and that 
$\ucI' = (M',N',E',\ucI',d')$ is a marked ideal such that $\ucI'$ is
equivalent to $\vp^*(\ucI)$ (in particular, $E' = \vp^{-1}(E)$) and
$\dim N' = \dim \vp^{-1}(N)$. Then $\vp$ lifts to smooth morphisms 
$\vp_j$ throughout the resolution towers for $\ucI'$ and $\ucI$, and,
for each $j$, $\inv_{\ucI'}=\inv_{\ucI}\circ\vp_j$, 
$\mu_{\ucI'}=\mu_{\ucI}\circ\vp_j$ and $J_{\ucI'}=J_{\ucI}\circ\vp_j$.
\end{enumerate}
\end{theorem}

\begin{proof}
We follow the steps of the proof of Theorem 1.3 in Section 5.
\smallskip

\noindent
{\bf Step I. $\ucI$ of maximal order.} 
\smallskip

\noindent
{\bf Case A. $E = \emptyset$.}
The coefficient ideals $\ucC_U(\ucI)$ admit resolution sequences and
invariants as required, by induction. Over each $U$, the desingularization
of $\ucI$ is realized by that of $\ucC_U(\ucI)$. If $x$ is a point of
$\cosupp \ucI_j$ lying over $U$, we set
$$
\oinv_{\ucI}(x) := \left(0, \inv_{\ucC_U(\ucI)}(x)\right),\quad
\omu_{\ucI}(x) := \mu_{\ucC_U(\ucI)}(x),\quad 
\oJ_{\ucI}(x) := J_{\ucC_U(\ucI)}(x).
$$

\noindent
{\bf Case B. General maximal order case.} The marked ideal $\ucJ$ (defined
using $\ucC_U(\ucI_{\emptyset})$) has empty exceptional divisor. We set
$$
\oinv_{\ucI}(x) := \left(s, \inv_{\ucJ}(x)\right),\quad
\omu_{\ucI}(x) := \mu_{\ucJ}(x),\quad
\oJ_{\ucI}(x) := J_{\ucJ}(x),
$$
for $x \in \cosupp \ucJ_j$, until $\ucJ$ is resolved. If the sequences of
blowings-up of the $\ucJ$ do not resolve $\ucI$, then we repeat the process
using the new $s = \max\,s(x)$, as in Section 5. 

Note that if $x \in \cosupp \ucI_j$ and $s(x) < s$, then all blowings-up in
the desingularization tower of $\ucI$ are isomorphisms over $x$, until we
reach a year $q_k$ where $s(x)$ is the maximum value. Then the values of
the invariants at $x$ equal their values in year $q_k$ (by definition;
see also Step II.B below).
\smallskip

\noindent
{\bf Step II. General $\ucI$.} First suppose that $\cI = 0$. If 
$x \in \cosupp \ucI = N$,
then we set
$$
\inv_{\ucI}(x) := \infty,\quad \mu_{\ucI}(x) := \infty,\quad 
J_{\ucI}(x) := \emptyset.
$$
We blow up with centre $N$ to resolve singularities.

Now suppose that $\cI \neq 0$.
\smallskip

\noindent
{\bf Case A. Monomial case.} If $x \in \cosupp \ucI_j$, then we set
$$
\inv_{\ucI}(x) := 0,\quad \mu_{\ucI}(x) := \mu_x(\ucI_j) = 
\frac{\ord_x\cI_j}{d},
$$
and we let $J_{\ucI}(x)$ denote the maximum among the subsets of $E_j$
that define the components of $\cosupp \ucI_j$ at $x$. The required properties
in this case have been proved in Section 4.
\smallskip

\noindent
{\bf Case B. General case.} We resolve the singularities of the 
companion ideal $\ucG := \ucG(\ucI)$ (which is of maximal order)
in order to reduce the maximum
order of the ideal $\cN(\ucI)$. If $x \in \cosupp \ucG_j$, then we set
$$
\inv_{\ucI}(x) := \left(\frac{\ord\,\cN(\ucI)}{d}, \oinv_{\ucG}(x)\right),\quad
\mu_{\ucI}(x) := \omu_{\ucG}(x),\quad
J_{\ucI}(x) := \oJ_{\ucG}(x).
$$

If $x \in \cosupp \ucI \setminus \cosupp \ucG$, then $\ord_x \cN(\ucI)$
will be the maximum order of $\cN(\ucI)$ in some neighbourhood of $x$,
so we can define the invariant in the same way over such a neighbourhood.
(When we resolve the singularities of $\ucG(\ucI)$, either 
$\cosupp \ucI_{r_1} = \emptyset$ or $\ord\,\cN(\ucI_{r_1}) <
\ord\,\cN(\ucI)$ after $r_1$ steps, and we repeat the process 
using $\left(\cN(\ucI_{r_1}),
\ord\,\cN(\ucI_{r_1})\right)$. If $x \in \cosupp \ucG(\ucI_{r_1})$ maps to
the complement of $\cosupp \ucG(\ucI)$, then all previous blowings-up are
isomorphisms at the successive images of $x$, so the values of the invariants
at these points are the same as at $x$.)

All properties required in the theorem follow by induction and
semicontinuity of $\ord_x$.
\end{proof}

\section{Comparison of algorithms and an example}
In this final section, we prove Corollary 1.4 and we give an example
to show that the versions of canonical desingularization (Theorem 1.1)
in \cite{EVnew, W, Ko} involve blowings-up with centres that are not
necessarily smooth.

\subsection{W{\l}odarczyk's method} Let $\ucI = (M,N,E,\cI,d)$ denote
a marked ideal. The objects constructed in Sections 3 and 4 using our
logarithmic derivative ideals $\ucD_E^j(\ucI)$ of course have analogues
defined using the standard derivative ideals $\ucD^k(\ucI)$ (cf. \cite{V2}):
We define $\ucC^k(\ucI) := \sum_{j=0}^k \ucD^j(\ucI)$, $k \leq d-1$. 
Then $\ucI$ and $\ucC^k(\ucI)$ are weak-equivalent (cf. Cor. 3.11),
but in general they are not equivalent \cite[Ex. 5.14]{BMda1}. For $\ucI$
of maximal order, W{\l}odarczyk
introduces a \emph{homogenized ideal} $\ucH(\ucI) = (M,N,E,\cH(\ucI),d)$,
where 
$$
\ucH(\ucI) := \sum_{j=0}^{d-1} \cD^j(\cI)\cdot(\cD^{d-1}(\cI))^j.
$$
Then $\ucH(\ucI)$ is weak-equivalent to $\ucI$ \cite[Lemma 3.5.2]{W}.
In general, $\ucH(\ucI)$ is not equivalent to $\ucI$. It follows from
Corollary 3.6, however, that if $E = \emptyset$, then $\ucC^k(\ucI)$
and $\ucH(\ucI)$ are both equivalent to $\ucI$. (The same statements hold
for the variant $\ucW(\ucI)$ of $\ucH(\ucI)$ introduced by Kollar \cite{Ko}.)

The desingularization algorithm as described in Section 5 allows many
minor variations. (See \cite{BMda1}.) For example, Step I.B (moving apart
$E$ and $\cosupp \ucI$) is described in a different way in \cite{W} that
may mean separating the components of $E$ from $\cosupp \ucI$ in a somewhat
different order. We will ignore this minor difference. The essential difference
in W{\l}odarczyk's treatment is his use of the homogenized ideal. Instead
of using the ideal $\ucC_E^{d-1}(\ucI)$ of \S3.4, he uses $\ucC^{d-1}(\ucH(\ucI))$.
Following the algorithm as presented in Section 5, this introduces a change
in Step I, where, instead of our coefficient ideal 
$\ucC_{\emptyset,P}(\ucI)$ (see (4.1)), W{\l}odarczyk would use
$$
\ucC_P(\ucI) := \left(U,P,\emptyset, \cC^{d-1}(\ucH(\ucI))|_P, d_{\ucC}\right).
$$
(Likewise in Kollar's version, using $\ucW(\ucI)$.)

\emph{A priori}, this change might result in a change in the actual recipe
for choosing the centres of blowing up. However, since $\ucC^{d-1}(\ucH(\ucI))$
is equivalent to $\ucI$ when $E = \emptyset$, it follows that $\ucC_P(\ucI)$
is equivalent to $\ucC_{\emptyset,P}(\ucI)$. We thus obtain Corollary 1.4.

\subsection{Example} Let $X$ denote the subvariety of $\IA^4$ defined by the
ideal $\cI$ generated by $y^2 - x^3$ and $x^4 + xz^2 - w^3$. We consider
resolution of singularities (Section 5) of the corresponding marked ideal $\ucI =
(\IA^4, \IA^4, \emptyset, \cI, 1)$ (which we write simply as 
$(\emptyset, \cI, 1)$). The ideal $\cI$ has maximum order $2$ precisely at
the origin, so that $C_0 := \{0\}$ is the centre of the first blowing-up
$\s_1$. Let $\ucI_1 = (E_1, \cI_1, 1)$ denote the transform of $\ucI$ by $\s_1$.
In the $x$-chart of $\s_1$, 
$$
\cI_1 = \left(x - y^2,\ x(x + z^2 - w^3)\right)
      = \left(x - y^2,\ y^2(y^2 + z^2 - w^3)\right),
$$ 
and $E_1$ is given by $(x)$. There is no nontrivial monomial part to factor
out, and $\ucI_1$ is of maximal order. According to Step I.B, the centre $C_1$
of the next blowing-up $\s_2$ is given by $x = y = 0$.

But the ideal of the strict transform $X_1$ of $X$ in the $x$-chart above
is generated by
$$
x - y^2,\quad y^2 + z^2 - w^3.
$$
Therefore, $C_1 \cap X_1$ is given by
$$
x = y = 0,\quad z^2 - w^3 = 0.
$$
Moreover, $\Sing X_1 = \{0\}$.

\bibliographystyle{ams}

\end{document}